\UseRawInputEncoding 
\documentclass[11pt, a4paper]{article}
\usepackage{amsfonts}
\usepackage{amssymb}
\usepackage{mathrsfs}
\usepackage{amsmath}
\usepackage{booktabs}
\usepackage{epsf, epsfig, amsfonts, amsgen, indentfirst}
\usepackage{amsmath, amstext, amsbsy, amsopn, amsthm, bbding, wasysym}
\usepackage{multicol, mathdots}
\usepackage{setspace}
\usepackage[numbers,sort&compress]{natbib}

\allowdisplaybreaks

\makeatletter
\renewcommand{\@makefntext}[1]{\noindent #1}
\makeatother

\begin{document}
\theoremstyle{plain}
\newtheorem{theorem}{ Theorem\ }[section]
\newtheorem{corollary}{Corollary\ }[section]
\newtheorem{lemma}{Lemma\ }[section]
\newtheorem{prop}{Proposition\ }[section]
\newtheorem{defi}{Definition\ }[section]
\newtheorem{property}{Property\ }[section]
\newtheorem{remark}{Remark\ }[section]
\newtheorem{example}{Example\ }[section]
\newtheorem{exercise}{Exercise\ }[section]
\newtheorem{conjecture}{Conjecture\ }[section]
\newtheorem{problem}{Problem\ }[section]

\renewcommand{\proofname}{\bf Proof}

\renewcommand{\proofname}{Proof}

\newenvironment{proofthm}[1]
{\begin{proof}[\textbf{#1}]}
{\end{proof}}
\renewcommand{\figurename}{Figure}
\renewcommand\refname{\centering References}

\allowdisplaybreaks

\baselineskip 15pt
\renewcommand{\baselinestretch}{1.2}
\title
{\LARGE \textbf{A sharp upper bound on the number of spanning forests of regular graphs
\thanks{ supported by NSFC (No. 12261071) and NSF of Qinghai Province (No. 2025-ZJ-902T).} 
}}

\author{Tingzeng Wu$^{a,b}$\thanks{{Corresponding author.\newline
\emph{E-mail address}: mathtzwu@163.com, sxlu863@163.com, a3566293588@163.com
}}, Shuxin Lu$^a$, Xiangshuai Dong$^c$ \\
{\small $^{a}$ School of Mathematics and Statistics, Qinghai Minzu University, }\\
{\small  Xining, Qinghai 810007, P.R.~China} \\
{\small $^{b}$ Qinghai Institute of Applied Mathematics,   Xining, Qinghai, 810007, P.R.~China}\\
{\small $^{c}$ School of Mathematical Sciences, Xiamen University,  }\\
{\small  Fujian, Xiamen 361005,  P.R.~China} \\
}

\date{}
\maketitle

\noindent {\bf Abstract: } 
Let $G$ be a simple graph on $n$ vertices, and let $F(G)$ denote the number of its spanning forests. 
Bencs and Csikvári [Upper bound for the number of spanning forests of regular graphs, European J. Combin. 110 (2023) 103677] proved that every $r$-regular graph $G$ with $r\geq 2$ satisfies $F(G) \leq r^{n}$. 
They further conjectured that for $r \geq 3$, 
\[
F(G)^{1/n} \leq \frac{(r - 1)^{r-1}}{(r^2 - 2r - 1)^{r/2-1}}.
\]
In this paper, we resolve this conjecture in the affirmative.
\\



\section{Introduction}

Throughout this paper, we consider only finite, simple, undirected graphs.
Let $G = (V(G), \\E(G))$ be a graph with vertex set $V(G) = \{v_1, v_2, \dots, v_n\}$
and edge set $E(G) = \{e_1, e_2, \dots, e_m\}$.
Let $|V(G)| = n$ and $|E(G)| = m$ denote the number of vertices and edges, respectively.
The degree of $v_i$ is the number of edges in $G$ incident to $v_i$, denoted by $d_i$.    
The sequence $(d_1, d_2, \dots, d_n)$ denotes the degree sequence of $G$.
Let $\Delta$ and $\delta$ denote the maximum and minimum degrees of $G$, respectively.
If $\Delta = \delta = r$, then $G$ is called an $r$-regular graph.
For notation not defined here, we refer the reader to \cite{Diestel2017}.

Let $S=\{v_1, \ldots, v_\ell\}$ be a nonempty subset of $V(G)$.
The subgraph induced by $S$, denoted by $G[S]$, has vertex set $S$ and edge set consisting of all edges of $G$ with both endpoints in $S$.
For $1 \leq i \leq \ell$, let $G_i = G - \{v_1, \ldots, v_{i-1}\}$; in particular, $G_1 = G$.
For any graph $N$ with $V(N) \subseteq V(G)$, the graph $G - V(N)$ is the subgraph induced by $V(G) \setminus V(N)$. For a single vertex $v_i$, we write $G - v_i$ for $G - \{v_i\}$.

A spanning forest of a graph \(G\) is a spanning subgraph \((V(G),E')\) that is acyclic, where \(E'\subseteq E(G)\). Write
\[
\mathcal{F}(G)=
\{(V(G),E'):E'\subseteq E(G),\ (V(G),E')\text{ is acyclic}\},
\qquad F(G)=|\mathcal{F}(G)|.
\]
Thus $F(G)$ denotes the number of spanning forests of $G$.

Counting spanning forests is a classical problem in combinatorial enumeration,
closely related to several areas of research, including random networks, probability theory, statistical physics, and network science \cite{Aldous2007, Annan1994, Borbényi2021, Caracciolo2004, Grimmett2004, Liu1981, Teranishi2005, ke2026}.
Jaeger et al. \cite{Jaeger1990} proved that computing $F(G)$ is \#P-hard.
In 2023, Bencs and Csikvári studied the number of spanning forests of regular graphs and established the following upper bound.
\begin{theorem}[\cite{Bencs2023}]
Let $G$ be an $r$-regular graph on $n$ vertices $(r \geq 2)$. Then
\[
F(G) \leq r^{n}.
\]
\end{theorem}
Based on Theorem 1.1, they proposed the following conjecture.
\begin{conjecture}
Let $G$ be an $r$-regular graph on $n$ vertices $(r \geq 3)$. Then
\[
F(G)^{1/n} \leq \frac{(r - 1)^{r-1}}{(r^2 - 2r - 1)^{r/2-1}}.
\]
\end{conjecture}
In this paper, we prove this conjecture, which is stated in the following theorem.
\begin{theorem}\label{thm2}
Let $G$ be an $r$-regular graph on $n$ vertices $(r \geq 3)$. Then
\[
F(G)^{1/n} \leq \frac{(r - 1)^{r-1}}{(r^2 - 2r - 1)^{r/2-1}}.
\]
\end{theorem}
\begin{remark}[Optimality of the exponential constant]\label{rem:sharpness}
Fix an integer $r\geq3$, and let $\mathcal{G}_r$ be the family of all nonempty finite simple $r$-regular graphs.
Choose a sequence $(G_k)_{k\geq1}$ in $\mathcal{G}_r$ whose girth, that is, the length of a shortest cycle, tends to infinity.
Such sequences exist; see the discussion of graph covers in \cite{Bencs2023}.
The limiting formula of Bencs and Csikv\'ari \cite{BencsTutte2022} gives
\[
\lim_{k\to\infty}F(G_k)^{1/|V(G_k)|}
=\frac{(r-1)^{r-1}}{(r^2-2r-1)^{r/2-1}}.
\]
Combining this limit with Theorem \ref{thm2} yields
\[
\sup_{G\in\mathcal{G}_r}F(G)^{1/|V(G)|}
=\frac{(r-1)^{r-1}}{(r^2-2r-1)^{r/2-1}}.
\]
Thus the exponential constant in Theorem \ref{thm2} is best possible.
As shown in Section 3, the inequality is nevertheless strict for every finite nonempty graph in $\mathcal{G}_r$.
\end{remark}
The paper is organized as follows. Section 2 presents the necessary lemmas and definitions. Section 3 provides the proof of Theorem 1.2. Section 4 concludes the paper with a brief summary.

\section{Preliminaries}

In this section, we introduce the definitions, notation, and lemmas needed for the proof of Theorem 1.2.

A matching of $G$ is a subset $M \subseteq E(G)$ in which no two edges are adjacent.
Let $|M|$ denote the number of edges in $M$, and let $V(M)$ denote the set of vertices covered by $M$.
Let $\mathcal{C}(G)$ denote the set of all cycles in $G$, and for $C \in \mathcal{C}(G)$, let $|C|$ denote its length.
Let $\mathcal{Q} \subseteq \mathcal{C}(G)$ be a set of pairwise vertex-disjoint cycles, and let $|\mathcal{Q}|$ denote the number of cycles in $\mathcal{Q}$.
Let $\mathcal{M}(G)$ denote the set of all matchings of $G$, including the empty matching.
Let $M' \in \mathcal{M}(G)$ be a matching satisfying $V(M') \cap \bigcup_{C \in \mathcal{Q}} V(C) = \varnothing$.
Define
\begin{equation}\label{eq:P_G}
P_G(x,t) = \sum_{(M',\mathcal{Q})} (-2)^{|M'|+|\mathcal{Q}|} t^{|\mathcal{Q}|} x^{|V(G)|-2|M'|-\sum_{C\in\mathcal{Q}}|C|},
\end{equation}
where the sum is over all pairs $(M',\mathcal{Q})$ satisfying the above conditions.
We adopt the conventions $P_{\varnothing}(x,t)=1$ and $0^0=1$.
When $t=0$, \eqref{eq:P_G} yields a matching polynomial in which each edge has weight $-2$, namely,
\[
M_G(x) = P_G(x,0) = \sum_{M \in \mathcal{M}(G)} (-2)^{|M|} x^{|V(G)|-2|M|}.
\]
For any vertex set $R \subseteq V(G)$, let $|R|$ denote its cardinality.
\begin{lemma}[\cite{Chaiken1982}]\label{lem:10}
Let $G$ be a graph on $n$ vertices, let $L(G)$ be its Laplacian matrix, and let $F$ be a spanning forest of $G$.
Let $L(R)$ denote the principal submatrix obtained from $L(G)$ by deleting the rows and columns indexed by the vertices in $R$, with the determinant of the empty principal submatrix defined to be $1$.
Then
\[
\det L(R) = \sum_{F} 1,
\]
where the sum is over all spanning forests $F$ of $G$ such that

\textup{(i)} $F$ has exactly $|R|$ connected components;

\textup{(ii)} each connected component $T$ of $F$ contains exactly one vertex of $R$.
\end{lemma}

\begin{lemma}\label{lem:11}
Let $G$ be a graph on $n$ vertices, and let $L(G)$ be its Laplacian matrix. Let
\[
Y=L(G)+\operatorname{diag}(c_1, c_2, \dots, c_n),
\]
where $c_i \geq 0$ for $1 \leq i \leq n$.
Then, for any vertex set $R\subseteq V(G)$, we have
\[
\det Y(R)
=
\sum_{\substack{F\in\mathcal{F}(G)\\
|V(T)\cap R|\le 1 \text{ for every component } T \text{ of } F}}
\prod_{\substack{T \text{ component of } F\\
V(T)\cap R=\varnothing}}
\left(\sum_{v_i \in V(T)}c_i\right),
\]
where $Y(R)$ denotes the principal submatrix obtained from $Y$ by deleting the rows and columns indexed by the vertex set $R$, and an empty product is defined to be $1$.
\end{lemma}
\begin{proof}
Let $S=V(G)\setminus R$ and write $S=\{v_{i_1}, \ldots, v_{i_j}\}$ in the original vertex order, where $j=|S|$.
For convenience, let
\[
D_S=\operatorname{diag}(c_{i_j}: v_{i_j}\in S).
\]
Thus
\[
Y(R)=L(R)+D_S.
\]
For $1\leq a\leq j$ and each $v_{i_a}\in S$, the column of $L(R)+D_S$ corresponding to $v_{i_a}$ is
\[
L_{i_a}+c_{i_a}e_{i_a},
\]
where $L_{i_a}$ denotes the $a$-th column of $L(R)$, and $e_{i_a}$ denotes the column vector whose $a$-th entry is $1$ and whose remaining entries are all $0$.
In expanding $\det(L(R)+D_S)$, there are two choices for each $v_{i_a}\in S$: choose $L_{i_a}$ or $c_{i_a}e_{i_a}$.

Fix a subset $U\subseteq S$, and choose $c_{i_a}e_{i_a}$ for vertices $v_{i_a}$ in $U$ and $L_{i_a}$ for vertices in $S\setminus U$.
The term arising from these choices contains the factor $\prod_{v_{i_a}\in U}c_{i_a}$.
Let $J_U=\{a: v_{i_a}\in U\}$.
With $U\subseteq S$ fixed, perform a Laplace expansion along the rows and columns corresponding to $U$.
Since both the deleted row set and the deleted column set are $U$, the cofactor sign in the Laplace expansion is
\[
(-1)^{\sum_{a\in J_U}a+\sum_{a\in J_U}a}=1.
\]
The remaining principal minor is $L(R\cup U)$.
Therefore
\begin{equation}\label{eq:20}
\det(L(R)+D_S)
=
\sum_{U\subseteq S}
\left(
\prod_{v_i\in U}c_i
\right)
\det L(R\cup U).
\end{equation}

By Lemma \ref{lem:10},
$\det L(R\cup U)$
equals the number of spanning forests $F$ satisfying the following conditions:

(i) $F$ has exactly $|R|+|U|$ connected components;

(ii) each connected component $T$ contains exactly one vertex of $R\cup U$, that is,
\[
|V(T)\cap(R\cup U)|=1.
\]
As seen from \eqref{eq:20}, the sum corresponds to first selecting $U\subseteq V(G)\setminus R$,
then selecting a spanning forest with root set $R\cup U$ and exactly one root in each connected component,
and assigning this object the weight $\prod_{v_i\in U}c_i$.
Since $R\subseteq R\cup U$, every contributing forest must satisfy
\begin{equation}\label{eq:21}
|V(T)\cap R|\leq 1
\qquad
\text{for every component } T.
\end{equation}
Now fix a spanning forest $F\in \mathcal{F}(G)$ satisfying \eqref{eq:21},
and divide the connected components of $F$ into two classes.

The first class consists of components $T'$ with
\[
V(T')\cap R\neq\varnothing.
\]
In this case, the component $T'$ contains exactly one vertex of $R$.
This vertex is already the only root permitted in the component within the root set $R\cup U$,
so $U$ cannot select any additional vertex from this component. Thus $T'$ contributes no factor $c_i$.

The second class consists of components $T''$ with
\[
V(T'')\cap R=\varnothing.
\]
For $T''$ to contain exactly one vertex of the root set $R\cup U$,
the set $U$ must select exactly one vertex from $V(T'')$. If the selected vertex is
$v_i\in V(T'')$, its contribution has weight $c_i$.
Thus, summing over all possible choices for this fixed component $T''$ gives the total weight
\begin{equation}
\sum_{v_i\in V(T'')}c_i.
\end{equation}

Roots in distinct connected components disjoint from $R$ can be chosen independently.
For a fixed forest $F$, sum over all admissible sets $U$: each such set contains no vertex from a component intersecting $R$ and exactly one vertex from each component disjoint from $R$.
The resulting total weight is

\begin{equation}\label{eq:23}
\prod_{\substack{
T \text{ component of } F\\
V(T)\cap R=\varnothing
}}
\left(
\sum_{v_i\in V(T)}c_i
\right).
\end{equation}

Selecting one root vertex in each component disjoint from $R$ determines the set $U$ uniquely.
This gives a weight-preserving bijection between forests satisfying \eqref{eq:21} with these selected roots and the pairs $(U,F)$ contributing to \eqref{eq:20}.
Combining \eqref{eq:20} and \eqref{eq:23} gives
\[
\det Y(R)
=
\sum_{\substack{F\in\mathcal{F}(G)\\
|V(T)\cap R|\le 1 \text{ for every component } T \text{ of } F}}
\prod_{\substack{T \text{ component of } F\\
V(T)\cap R=\varnothing}}
\left(\sum_{v_i \in V(T)}c_i\right).
\]
This proves the lemma.
\end{proof}

In fact, Lemma \ref{lem:11} is a generalization of the matrix-forest theorem.
Setting $c_i=1$ for $1\leq i\leq n$
yields the classical matrix-forest theorem and its related forms; see \cite{Chebotarev1997}.
\begin{lemma}[\cite{Harary1962, Sachs1964}]
Let $G$ be a graph on $n$ vertices, let $A(G)$ be its adjacency matrix, and let $I$ be the identity matrix of order $n$.
Then
\begin{equation}\label{eq:4}
\det(xI-A(G))
=
\sum_{(M, \mathcal{Q})}
(-1)^{|M|+|\mathcal{Q}|}2^{|\mathcal{Q}|}
x^{|V(G)|-2|M|-\sum_{C\in \mathcal{Q}}|C|}.
\end{equation}
\end{lemma}
\begin{lemma}\label{lem:forest-identity}
Let $G$ be a graph on $n$ vertices, and let $x\geq 1+\Delta$.
Then
\begin{equation}\label{eq:forest-identity}
 P_G(x,1)
 =\sum_{F\in\mathcal{F}(G)}
   \prod_{T \text{ component of } F}
   \left[1+\sum_{v_i\in V(T)}(x-1-d_i)\right].
\end{equation}
If $G$ is an $r$-regular graph, then
\begin{equation}\label{eq:F-exact}
 P_G(r+1,1)=F(G).
\end{equation}
If $H$ is any induced subgraph of an $r$-regular graph $G$ and $x\geq r+1$, then $P_H(x,1)>0$.
\end{lemma}

\begin{proof}
For a fixed $N\in\mathcal{M}(G)$,
applying \eqref{eq:4} to the graph $G-V(N)$ gives
\begin{align}
(-1)^{|N|}
&\det(xI-A(G-V(N))) \notag\\
&\quad
=\sum_{(M',\mathcal{Q})}
(-1)^{|N|+|M'|+|\mathcal{Q}|}2^{|\mathcal{Q}|}
x^{|V(G)|-2|N|-2|M'|-\sum_{C\in\mathcal{Q}}|C|}.
\label{eq:5}
\end{align}
Here $M'\in\mathcal{M}(G-V(N))$,
and $\mathcal{Q}$ is a set of pairwise vertex-disjoint cycles in $G-V(N)$ such that
\[
V(M')\cap\bigcup_{C\in\mathcal{Q}}V(C)=\varnothing.
\]

Let $\widetilde{M}=N \cup M'$.
Since $M'$ and $N$ are vertex-disjoint,
$\widetilde{M}$ is still a matching of $G$,
and $\widetilde{M}\in\mathcal{M}(G)$.
Conversely, for any fixed $(\widetilde{M},\mathcal{Q})$,
assign each edge of $\widetilde{M}$ independently to either $N$ or $M'$.
Each such choice uniquely determines a decomposition $\widetilde{M}=N \cup M'$.
There are therefore $2^{|\widetilde{M}|}$ such decompositions, and each has sign
$(-1)^{|N|+|M'|}=(-1)^{|\widetilde{M}|}$. Hence
\[
\sum_{N \cup M'=\widetilde{M}}
(-1)^{|N|+|M'|}
=2^{|\widetilde{M}|}(-1)^{|\widetilde{M}|}
=(-2)^{|\widetilde{M}|}.
\]
First summing \eqref{eq:5} over all $N\in\mathcal{M}(G)$,
and then combining all terms with the same $(\widetilde{M},\mathcal{Q})$,
gives
\begin{equation}\label{eq:11}
P_G(x,1)
=
\sum_{N\in\mathcal{M}(G)}
(-1)^{|N|}
\det(xI-A(G-V(N))).
\end{equation}

Let $c_i=x-d_i$, where $d_i$ is the degree of the vertex $v_i$ in $G$.
Thus $c_i \geq 1$ and $xI-A(G)=L(G)+\operatorname{diag}(c_1, c_2, \dots, c_n)$.

For a fixed $N\in\mathcal{M}(G)$, apply Lemma \ref{lem:11} with $R=V(N)$.
Since $xI-A(G-V(N))=(xI-A(G))(V(N))$,
the corresponding determinant in \eqref{eq:11} can be expressed as
a sum over forests $F\in\mathcal{F}(G)$ satisfying the following conditions:
each component of $F$ contains at most one vertex of $V(N)$;
a component $T$ disjoint from $V(N)$ contributes the weight $\sum_{v_i\in V(T)}c_i$.
Now add the edges of the matching $N$ back to this forest $F$.
Since each tree component of $F$ contains at most one vertex of $V(N)$ and $N$ is a matching,
the endpoints of each edge of $N$ must lie in two distinct tree components,
and each original tree component is incident with at most one edge of $N$. Thus, adding $N$ back creates no cycle,
and the result is again a spanning forest $F' \in\mathcal{F}(G)$.
Consequently, after $N$ is added back, each resulting tree component contains at most one edge of $N$.
Conversely, for each tree component $T$ of the resulting spanning forest $F'$,
its origin can equivalently be described by one of the following two mutually exclusive types of marking:

(i) Select a vertex $v_i$ in $T$ as a root vertex: this marking has weight $c_i$,
so the total weight of all vertex-root markings is $\sum_{v_i\in V(T)}c_i$;

(ii) Select an edge in $T$ as a root edge: this edge belongs to the matching $N$, and the factor $(-1)^{|N|}$ in~\eqref{eq:11} gives each root-edge marking a contribution of $-1$.
Thus the total weight of all root-edge markings is $-|E(T)|$.

\noindent Deleting the root edge corresponding to (ii) restores the tree component of $F'$ to the two tree components of $F$ that each contain
one endpoint of that root edge.
Thus the above correspondence between ``a matching $N$ together with a rooted forest from Lemma \ref{lem:11}'' and
``a spanning forest marked with vertex roots or root edges'' is a bijection.
For a fixed tree component $T$, the total weight of all markings is
\begin{align*}
\sum_{v_i\in V(T)}c_i-|E(T)|
=1+\sum_{v_i\in V(T)}(x-1-d_i).
\end{align*}

Since markings in distinct tree components are independent, summing over all spanning forests and taking the product of the component weights yields \eqref{eq:forest-identity}.

For an $r$-regular graph, set $x=r+1$; every vertex $v_i$ then satisfies $x-1-d_i=0$.
Thus the bracketed expression for each tree component in~\eqref{eq:forest-identity} equals $1$,
which gives~\eqref{eq:F-exact}.

Finally, if $H$ is an induced subgraph of an $r$-regular graph,
then $\deg_H(v_i)\leq r$ for every vertex $v_i \in V(H)$,
where $\deg_H(v_i)$ denotes the degree of $v_i$ in $H$.
When $x\geq r+1$,
the bracketed expression for each tree component in~\eqref{eq:forest-identity} is at least $1$, so $P_H(x,1)>0$.
\end{proof}

\begin{lemma}\label{lem:vertex-recurrence}
Let $P_G(x,t)$ be the polynomial defined in \eqref{eq:P_G}.
Then, for every $v_i\in V(G)$, the following recurrence holds:
\begin{align}\label{eq:vertex-recurrence}
 P_G(x,t)
 ={}&xP_{G-v_i}(x,t)
      -2\sum_{v_j\in N_G(v_i)}P_{G-\{v_j,v_i\}}(x,t)
    \notag\\
    &-2t\sum_{\substack{C\in\mathcal{C}(G)\\v_i\in V(C)}}
                   P_{G-V(C)}(x,t),
\end{align}
where $N_G(v_i)$ denotes the set of all vertices adjacent to $v_i$.
\end{lemma}
\begin{proof}
By the definition of $P_G(x,t)$, for a fixed vertex $v_i$,
exactly one of the following three cases occurs in any contributing object $(M',\mathcal{Q})$:

(i) $v_i\notin V(M')$ and $v_i\notin\bigcup_{C\in\mathcal{Q}}V(C)$;

(ii) $v_i\in V(M')$;

(iii) $v_i\in V(C)$, where $C\in \mathcal{Q}$.

Consider the first case. Suppose that $v_i$ belongs neither to an edge of the matching $M'$
nor to a cycle in $\mathcal{Q}$.
Then $(M',\mathcal{Q})$ lies entirely in the induced subgraph $G-v_i$.
Thus these objects are in bijection with the objects in $P_{G-v_i}(x,t)$.
If the weight of $(M',\mathcal{Q})$ in $P_{G-v_i}(x,t)$ is
\[
(-2)^{|M'|+|\mathcal{Q}|}t^{|\mathcal{Q}|}
 x^{|V(G-v_i)|-2|M'|-\sum_{C\in \mathcal{Q}}|C|},
\]
then, since $|V(G)|=|V(G-v_i)|+1$,
its weight in $P_G(x,t)$ differs from the above expression by exactly a factor of $x$.
Hence the total contribution of case (i) is $xP_{G-v_i}(x,t)$.

Now consider the second case.
Suppose that $v_i\in V(M')$. Since $M'$ is a matching, there is a unique vertex
$v_j\in N_G(v_i)$ such that $v_iv_j\in M'$. Let
$M''=M'\setminus\{v_iv_j\}$.
Since the edges of $M'$ are pairwise nonadjacent and $V(M')$ is disjoint from the vertices of all cycles in $\mathcal{Q}$,
after deleting $v_i,v_j$, the pair $(M'',\mathcal{Q})$ is precisely an object in the graph
$G-\{v_i,v_j\}$ that satisfies the required conditions.

Conversely, fix any $v_j\in N_G(v_i)$. If $(M'',\mathcal{Q})$ is
an admissible object in the graph $G-\{v_i,v_j\}$, then
$M'=M''\cup\{v_iv_j\}$
is still a matching of $G$, and the added edge $v_iv_j$ is vertex-disjoint from $M''$ and all cycles in $\mathcal{Q}$.
Thus, for a fixed $v_j\in N_G(v_i)$, the above correspondence is a bijection.
Under this correspondence,
\[
|M'|=|M''|+1,
\qquad
|V(G)|=|V(G-\{v_i,v_j\})|+2.
\]
The weight of the object $(M',\mathcal{Q})$ in $G$ is therefore
\begin{align*}
&(-2)^{|M'|+|\mathcal{Q}|}t^{|\mathcal{Q}|}
 x^{|V(G)|-2|M'|-\sum_{C\in \mathcal{Q}}|C|}
\\
&\quad=
(-2)^{|M''|+1+|\mathcal{Q}|}t^{|\mathcal{Q}|}
 x^{|V(G-\{v_i,v_j\})|-2|M''|-\sum_{C\in \mathcal{Q}}|C|}
\\
&\quad=
-2\,(-2)^{|M''|+|\mathcal{Q}|}t^{|\mathcal{Q}|}
 x^{|V(G-\{v_i,v_j\})|-2|M''|-\sum_{C\in \mathcal{Q}}|C|}.
\end{align*}
Hence, for a fixed neighbor $v_j$, the total contribution of these objects is
$-2P_{G-\{v_i,v_j\}}(x,t)$.
Summing over all $v_j\in N_G(v_i)$, the total contribution of case (ii) is
\[
-2\sum_{v_j\in N_G(v_i)}P_{G-\{v_i,v_j\}}(x,t).
\]

Finally, consider the third case.
Suppose that $v_i$ belongs to a cycle in $\mathcal{Q}$.
Since the cycles in $\mathcal{Q}$ are pairwise vertex-disjoint, the cycle containing $v_i$ is unique.
Denote this cycle by $C$. Let
$\mathcal{Q}'=\mathcal{Q}\setminus\{C\}$.
The vertex-disjointness of $\mathcal{Q}$, together with
\[
V(M')\cap\bigcup_{C'\in \mathcal{Q}}V(C')=\varnothing,
\]
implies that, after deleting $V(C)$, the pair $(M',\mathcal{Q}')$ is precisely an object in the graph
$G-V(C)$ satisfying the defining conditions.
Conversely, fix a cycle $C\in\mathcal{C}(G)$ satisfying $v_i\in V(C)$.
If $(M',\mathcal{Q}')$ is any admissible object in $G-V(C)$, then setting
$\mathcal{Q}=\mathcal{Q}'\cup\{C\}$
gives a unique object in $G$ belonging to case (iii). Thus,
for a fixed $C$, this correspondence is also a bijection.
Under this correspondence,
\[
|\mathcal{Q}|=|\mathcal{Q}'|+1,
\qquad
\sum_{C'\in \mathcal{Q}}|C'|
=|C|+\sum_{C'\in \mathcal{Q}'}|C'|,
\]
and $|V(G)|=|V(G-V(C))|+|C|$.
Hence the weight of $(M',\mathcal{Q})$ in $P_G(x,t)$ is
\begin{align*}
&(-2)^{|M'|+|\mathcal{Q}|}t^{|\mathcal{Q}|}
 x^{|V(G)|-2|M'|-\sum_{C'\in \mathcal{Q}}|C'|}
\\
&\quad=
(-2)^{|M'|+|\mathcal{Q}'|+1}t^{|\mathcal{Q}'|+1}
 x^{|V(G-V(C))|-2|M'|-\sum_{C'\in \mathcal{Q}'}|C'|}
\\
&\quad=
-2t\,(-2)^{|M'|+|\mathcal{Q}'|}t^{|\mathcal{Q}'|}
 x^{|V(G-V(C))|-2|M'|-\sum_{C'\in \mathcal{Q}'}|C'|}.
\end{align*}
Therefore, for a fixed cycle $C$, its total contribution is
$-2t\, P_{G-V(C)}(x,t)$.
Summing over all cycles $C\in\mathcal{C}(G)$ satisfying $v_i\in V(C)$,
the total contribution of case (iii) is
\[
-2t\sum_{\substack{C\in\mathcal{C}(G)\\v_i\in V(C)}}P_{G-V(C)}(x,t).
\]

Adding the contributions of the three cases gives
\[
P_G(x,t)
=
 xP_{G-v_i}(x,t)
-2\sum_{v_j\in N_G(v_i)}P_{G-\{v_i,v_j\}}(x,t)
-2t\sum_{\substack{C\in\mathcal{C}(G)\\v_i\in V(C)}}P_{G-V(C)}(x,t).
\]
This completes the proof.
\end{proof}

\begin{lemma}\label{lem:cycle-derivative}
Let $P_G(x,t)$ be the polynomial defined in \eqref{eq:P_G}. Then
\begin{equation}\label{eq:cycle-derivative}
 \frac{\partial}{\partial t}P_G(x,t)
 =-2\sum_{C\in\mathcal{C}(G)}P_{G-V(C)}(x,t).
\end{equation}
\end{lemma}
\begin{proof}
Taking the partial derivative of \eqref{eq:P_G} with respect to the cycle parameter $t$ gives
\[
\begin{aligned}
\frac{\partial}{\partial t}P_G(x,t)
&=
\sum_{\substack{(M',\mathcal{Q})\\\mathcal{Q}\neq\varnothing}}
(-2)^{|M'|+|\mathcal{Q}|}
|\mathcal{Q}|t^{|\mathcal{Q}|-1}
x^{|V(G)|-2|M'|-\sum_{C\in \mathcal{Q}}|C|}                                      \\
&=
\sum_{\substack{(M',\mathcal{Q})\\\mathcal{Q}\neq\varnothing}}
\sum_{C\in \mathcal{Q}}
(-2)^{|M'|+|\mathcal{Q}|}
t^{|\mathcal{Q}|-1}
x^{|V(G)|-2|M'|-\sum_{C\in \mathcal{Q}}|C|}.
\end{aligned}
\]
The factor $|\mathcal{Q}|$ arising from differentiation can be interpreted as
selecting and fixing a cycle $C$ from the set $\mathcal{Q}$.

Now fix $C''\in\mathcal{C}(G)$, and let $\mathcal{Q}''=\mathcal{Q}\setminus\{C''\}$.
Since the cycles in $\mathcal{Q}$ are pairwise vertex-disjoint,
and the matching $M'$ is disjoint from the vertices of all cycles in $\mathcal{Q}$,
after deleting all vertices of $C''$, the pairs $(M',\mathcal{Q}'')$ on $G-V(C'')$ are in bijection with
the pairs $(M',\mathcal{Q})$ on $G$ satisfying $C''\in\mathcal{Q}$.
Under this correspondence,
\[
|\mathcal{Q}|=|\mathcal{Q}''|+1,
\qquad
\sum_{C\in \mathcal{Q}}|C|
=
|C''|+\sum_{C\in \mathcal{Q}''}|C|,
\]
and $|V(G-V(C''))|=|V(G)|-|C''|$.
The weight of the corresponding term can therefore be rewritten as
\[
\begin{aligned}
&(-2)^{|M'|+|\mathcal{Q}|}
t^{|\mathcal{Q}|-1}
x^{|V(G)|-2|M'|-\sum_{C\in \mathcal{Q}}|C|}                                      \\
&\qquad
=
-2(-2)^{|M'|+|\mathcal{Q}''|}
t^{|\mathcal{Q}''|}
x^{|V(G-V(C''))|-2|M'|-\sum_{C\in \mathcal{Q}''}|C|}.
\end{aligned}
\]
Summing over all pairs $(M',\mathcal{Q}'')$ on $G-V(C'')$ gives $-2P_{G-V(C'')}(x,t)$.
Finally, summing over all cycles $C\in\mathcal{C}(G)$ in $G$ yields \eqref{eq:cycle-derivative}.
\end{proof}

\begin{lemma}\label{lem:monotonicity}
Fix an integer $z\geq3$ and a real number $x\geq z+1$.
For every graph $G$ with $\Delta \leq z$ and every $0\leq t\leq1$, we have $P_G(x,t)>0$.
Moreover, for $v_i\in V(G)$, the function
\begin{equation}\label{eq:R-definition}
 R_G(v_i,t)=\frac{P_G(x,t)}{P_{G-v_i}(x,t)}
\end{equation}
is nonincreasing in $t$.
\end{lemma}

\begin{proof}
We first prove $P_G(x,t)>0$ by induction on the number of vertices $|V(G)|$.

For the empty graph, $P_\varnothing(x,t)=1$, so the assertion holds.
Suppose that the assertion holds for all proper induced subgraphs.
Then \eqref{eq:cycle-derivative} gives $\partial P_G(x,t)/\partial t\leq0$.
By Lemma \ref{lem:forest-identity}, we have $P_G(x,1)>0$,
so $P_G(x,t)\geq P_G(x,1)>0$, completing the proof.
Next, we again use induction on the number of vertices $|V(G)|$ to prove that $R_G(v_i,t)$ is nonincreasing.

When $|V(G)|=1$, let $V(G)=\{v_1\}$. There are no matching edges or cycles, so
\[
P_G(x,t)=x, \qquad P_{G-v_1}(x,t)=P_\varnothing(x,t)=1.
\]
$R_G(v_1,t)=x$,
which is independent of $t$ and is clearly nonincreasing.

Now suppose that the assertion holds for all proper induced subgraphs.
Let $K$ be a proper induced subgraph of $G$, and let $S=\{s_1,\ldots,s_k\}\subseteq V(K)$.
By \eqref{eq:R-definition}, we obtain the following chain factorization:
\begin{equation}\label{eq:ratio-chain}
 \frac{P_{K-S}(x,t)}{P_K(x,t)}
 =\prod_{i=1}^{k}
 \frac{1}{R_{K-\{s_1,\ldots,s_{i-1}\}}(s_i,t)}.
\end{equation}
For $i=1$, we adopt the convention
$R_{K-\{s_1,\ldots,s_{i-1}\}}(s_i,t)=R_K(s_1,t)$.
By the induction hypothesis, for every $i$ $(1\leq i \leq k)$,
$R_{K-\{s_1,\ldots,s_{i-1}\}}(s_i,t)>0$, and this ratio is nonincreasing in $t$.
Thus $\frac{P_{K-S}(x,t)}{P_K(x,t)}$ is nondecreasing in $t$.

Dividing \eqref{eq:vertex-recurrence} by $P_{G-v_i}(x,t)$ gives
\begin{equation}\label{eq:6}
\begin{split}
\frac{P_G(x,t)}{P_{G-v_i}(x,t)} &= x - 2 \sum_{v_j \in N_G(v_i)} \frac{P_{G-\{v_i,v_j\}}(x,t)}{P_{G-v_i}(x,t)} \\
&\quad - 2t \sum_{\substack{C \in \mathcal{C}(G) \\ v_i \in V(C)}} \frac{P_{G-V(C)}(x,t)}{P_{G-v_i}(x,t)}.
\end{split}
\end{equation}
Since $\frac{P_{K-S}(x,t)}{P_K(x,t)}$ is nondecreasing in $t$,
both $\frac{P_{G-\{v_i,v_j\}}(x,t)}{P_{G-v_i}(x,t)}$ 
and $\frac{P_{G-V(C)}(x,t)}{P_{G-v_i}(x,t)}$ are nondecreasing in $t$.
Therefore $R_G(v_i,t)$ is nonincreasing in $t$, completing the proof.
\end{proof}
The monotonicity of $R_G(v_i,t)$, together with \eqref{eq:ratio-chain}, gives the following corollary.

\begin{corollary}\label{cor:1}
Fix an integer $z\geq3$ and a real number $x\geq z+1$.
Let $G$ be a graph satisfying $\Delta \leq z$, and let $0\leq t\leq1$.
Then, for every $S\subseteq V(G)$, the function $P_{G-S}(x,t)/P_G(x,t)$ is nondecreasing in $t$.
\end{corollary}

\begin{lemma}\label{lemma:cycle-loss}
Let $G$ be an $r$-regular graph on $n$ vertices $(r\geq3)$. Let $x_0=r+1$ and
$\rho_r=x_0-r/x_0$.
Then
\begin{equation}\label{eq:cycle-loss}
 \ln F(G)
 \leq\ln M_G(x_0)-2\sum_{C\in\mathcal{C}(G)}\rho_r^{-|C|}.
\end{equation}
\end{lemma}
\begin{proof}
The rules for partial differentiation give
$\frac{\partial}{\partial t} \ln P_G(x_0,t) = \frac{1}{P_G(x_0,t)} \frac{\partial P_G(x_0,t)}{\partial t}$,
and combining \eqref{eq:cycle-derivative} with Corollary \ref{cor:1} gives
\begin{equation}\label{eq:loss-first}
\ln M_G(x_0) - \ln F(G) 
= 2 \sum_{C \in \mathcal{C}(G)} \int_0^1 
\frac{P_{G-V(C)}(x_0,t)}{P_G(x_0,t)} dt  
\ge 2 \sum_{C \in \mathcal{C}(G)} 
\frac{M_{G-V(C)}(x_0)}{M_G(x_0)}.
\end{equation}
It therefore suffices to prove
\begin{equation}\label{eq:25}
\sum_{C \in \mathcal{C}(G)} 
\frac{M_{G-V(C)}(x_0)}{M_G(x_0)}
\geq \sum_{C\in\mathcal{C}(G)}\rho_r^{-|C|}.
\end{equation}
For any nonempty induced subgraph $H$ of $G$,
let $\rho_H(v_i)=\frac{M_H(x_0)}{M_{H-v_i}(x_0)}$.
Combining Lemma \ref{lem:monotonicity} with \eqref{eq:vertex-recurrence} at $t=0$ gives
\[
 0<\rho_H(v_i)
 =x_0 - 2 \sum_{v_j\in N_H(v_i)} \frac{M_{H-\{v_j,v_i\}}(x_0)}{M_{H-v_i}(x_0)}
 =x_0-2\sum_{v_j\in N_H(v_i)}\frac1{\rho_{H-v_i}(v_j)}\leq x_0.
\]
Hence $\frac{1}{\rho_{H-v_i}(v_j)} \geq \frac{1}{x_0}$,
and thus
\begin{equation}\label{eq:rho-upper}
 \rho_H(v_i) \leq x_0 - 2 \sum_{v_j\in N_H(v_i)}\frac1{x_0} = x_0-\frac{2\deg_H(v_i)}{x_0},
\end{equation}
where $\deg_H(v_i)$ denotes the degree of the vertex $v_i$ in $H$.
By \eqref{eq:rho-upper} and the arithmetic--geometric mean inequality,
\begin{equation}
\frac{M_G(x_0)}{M_{G-S}(x_0)}
=
\prod_{i=1}^{\ell}\rho_{G_i}(v_i)
\leq
\prod_{i=1}^{\ell}
\left(
x_0-\frac{2\deg_{G_i}(v_i)}{x_0}
\right)
\leq
\left(
x_0-\frac{2}{\ell x_0}
\sum_{i=1}^{\ell}\deg_{G_i}(v_i)
\right)^{\ell}.
\label{eq:minor-product}
\end{equation}
Each edge with both endpoints in $S$ is counted exactly once fewer during successive vertex deletions, so
\[
 \sum_{i=1}^{\ell}\deg_{G_i}(v_i)
 =r\ell-|E(G[S])|\geq\frac{r\ell}{2}.
\]
Substituting this into \eqref{eq:minor-product} gives
\begin{equation}\label{eq:minor-lower}
 \frac{M_{G-S}(x_0)}{M_G(x_0)}\geq \rho_r^{-\ell}.
\end{equation}
For any $C\in\mathcal{C}(G)$,
setting $S=V(C)$, $\ell=|C|$ in \eqref{eq:minor-lower} yields
\[
\frac{M_{G-V(C)}(x_0)}
{M_G(x_0)}
\ge
\rho_r^{-|C|}.
\]
Summing over all $C\in\mathcal{C}(G)$ gives \eqref{eq:25}.
Combining this with \eqref{eq:loss-first} proves Lemma \ref{lemma:cycle-loss}.
\end{proof}

Fix an integer $z\geq3$. If
\begin{equation}\label{eq:spectral-domain}
 x>2\sqrt{2(z-1)},
\end{equation}
then define the following functions.
\begin{equation}\label{eq:yqg}
 y=y(x)=\frac{x+\sqrt{x^2-8(z-1)}}2,
 \qquad q=q(x)=\frac2{y^2}, \qquad
 h=\frac1y, \qquad g=\frac{y}{y^2-2}.
\end{equation}
Thus
\begin{equation}\label{eq:y-identities}
 x=y+\frac{2(z-1)}y, \qquad
 0<q<\frac1{z-1}, \qquad
 h=\frac1{x-2(z-1)h}, \qquad
 g=\frac1{x-2zh}.
\end{equation}
Let
\begin{equation}\label{eq:tree-expression}
 \mathcal{A}_z(x)=\frac{y(x)^{z-1}}
                {(y(x)^2-2)^{z/2-1}}.
\end{equation}

\begin{lemma}[\cite{Heilmann1972}]\label{lem:30}
Let $G$ be a graph on $n$ vertices with matching polynomial
\[
\mu(G,x)
=
\sum_{k=0}^{\lfloor n/2\rfloor}
(-1)^k m_k(G)x^{n-2k},
\]
where $m_k(G)$ denotes the number of matchings of $G$ containing $k$ edges.
Then all zeros of $\mu(G,x)$ are real. Moreover, when $\Delta\ge 2$,
every zero $\theta$ of $\mu(G,x)$ satisfies
$|\theta|\le 2\sqrt{\Delta-1}$.
\end{lemma}

\begin{lemma}\label{lem:matching-positive}
When \eqref{eq:spectral-domain} holds, every graph $G$ with $\Delta \leq z$
satisfies $M_G(x)>0$.
If $v_i\in V(G)$ and $d_i\leq z-1$, then
\begin{equation}\label{eq:cavity-lower}
 \frac{M_G(x)}{M_{G-v_i}(x)}\geq y.
\end{equation}
\end{lemma}

\begin{proof}
Since $M_G(x)=2^{n/2}\mu(G,x/\sqrt2)$,
Lemma \ref{lem:30} implies that $M_G(x)>0$ for every graph $G$ satisfying $2 \leq \Delta \leq z$.
When $\Delta \le1$, the graph $G$ is a disjoint union of isolated vertices and independent edges, so
\[
M_G(x)=x^{|V(G)|-2|E(G)|}(x^2-2)^{|E(G)|}>0.
\]
Consequently, when \eqref{eq:spectral-domain} holds, every graph $G$ with $\Delta \leq z$
satisfies $M_G(x)>0$.
We now prove \eqref{eq:cavity-lower} by induction on the number of vertices $|V(G)|$.

When $|V(G)|=1$, let $V(G)=\{v_1\}$. There are no matching edges, so
$M_G(x)=x$, $M_\varnothing(x)=1$.
Thus
\[
\frac{M_G(x)}{M_\varnothing(x)}=x.
\]
Since $x=y+\frac{2(z-1)}y$, we have $x>y$, and \eqref{eq:cavity-lower} holds.

Suppose that \eqref{eq:cavity-lower} holds for all proper induced subgraphs of $G$.
For a nonempty graph $G$, setting $t=0$ in \eqref{eq:6} gives
\begin{equation}\label{eq:7}
 \frac{M_G(x)}{M_{G-v_i}(x)}
 =x-2\sum_{v_j\in N_G(v_i)}
       \frac{M_{G-\{v_j,v_i\}}(x)}{M_{G-v_i}(x)}.
\end{equation}
Since $v_j\in N_G(v_i)$,
we have $\deg_{G-v_i}(v_j)\leq z-1$ (where $\deg_{G-v_i}(v_j)$ denotes
the degree of $v_j$ in $G-v_i$);
the induction hypothesis implies that every summand in \eqref{eq:7} is at most $1/y$.
Therefore
\begin{equation}\label{eq:8}
 \frac{M_G(x)}{M_{G-v_i}(x)}
 \geq x-\frac{2d_i}y.
\end{equation}
If $d_i\leq z-1$, substituting this into \eqref{eq:8} proves Lemma \ref{lem:matching-positive}.
\end{proof}

\begin{defi}[\cite{Godsil1981}]
Given a vertex $v$ of $G$, the vertices of its path tree $\mathcal{T}_v$ are all sequences
$(v_0,v_1,\ldots,v_k)$, where $v_0=v$, the vertices are pairwise distinct, and consecutive vertices are adjacent in $G$.
The root is the singleton sequence $(v)$; the parent of a nonroot vertex is obtained by deleting the last term of its sequence.
\end{defi}

The above classical path-tree construction shows that $\mathcal{T}_v$ is a finite tree.
In what follows, let $G$ be an $r$-regular graph $(r\geq3)$,
set $z=r$ in \eqref{eq:yqg}--\eqref{eq:tree-expression},
and assume throughout that $x>2\sqrt{2(r-1)}$.
The root $(v)$ of the path tree $\mathcal{T}_v$ has $r$ children,
and every nonroot node has at most $r-1$ children.
For a vertex $\omega$ of the tree, define the following recursively from the leaves toward the root:
\begin{equation}\label{eq:h-recursion}
 h_\omega=\frac1{x-2\sum_{\omega'\text{ is }\omega\text{'s child}}h_{\omega'}}.
\end{equation}
For a leaf $\omega$, we have $h_\omega=1/x\leq h$.
If all values at the children of a nonroot node belong to $(0,h]$, then the denominator is at least
$x-2(r-1)h=y>0$, so the value at this node is also at most $h$.
The denominator at the root is at least $x-2rh=(y^2-2)/y>0$.
Thus all denominators in the above definition are positive,
and we denote the value at the root by $g_v$.
Godsil's path-tree identity \cite{Godsil1981} gives
\[
\mu(G-v,x)\mu(\mathcal{T}_v,x)
=
\mu(G,x)\mu(\mathcal{T}_v-(v),x).
\]
Applying the vertex-deletion recurrence to the matching polynomial of the path tree and proceeding inductively from the leaves to the root,
we obtain from \eqref{eq:h-recursion} that
\[
g_v=\frac{M_{\mathcal{T}_v-(v)}(x)}{M_{\mathcal{T}_v}(x)}.
\]
Then, substituting $x/\sqrt2$ for the variable $x$ in Godsil's path-tree identity and using
$M_G(x)=2^{n/2}\mu(G,x/\sqrt2)$,
gives
\begin{equation}\label{eq:root-ratio}
 g_v=\frac{M_{G-v}(x)}{M_G(x)}.
\end{equation}

Let a nonroot node of the path tree $\mathcal{T}_v$ be
$\omega=(v_0,v_1,\ldots,v_k)$,
where $k\ge 1, v_0=v$.
Define
$G_\omega
=
G-\{v_0,v_1,\ldots,v_{k-1}\}$.
By \eqref{eq:7}, we have
$h_\omega = \frac{M_{G_\omega - v_k}(x)}{M_{G_\omega}(x)}$.
Combining this with Lemma \ref{lem:matching-positive} and \eqref{eq:yqg} gives
\[
 0<h_\omega\leq h\quad(\omega\text{ is not the root}),
 0<g_v
=\frac1{x-2\sum_{\omega\text{ is a child of the root}}h_\omega}
\le \frac1{x-2rh}=g.
\]
For a nonroot vertex $\omega$, let $k_\omega$ be its number of children, and let
$a_\omega=r-1-k_\omega$.
$a_\omega$ denotes the number of edges from the endpoint of the current path to previously visited vertices, excluding an immediate return to the parent.

For the root $v$, let $\mathcal{B}_v$ be the set of vertex sequences of the form
\begin{equation}\label{eq:blocked-definition}
 \beta=(v_0,\ldots,v_{\ell-1},v_\ell),
 \qquad v_0=v,
\end{equation}
where $v_0,\ldots,v_{\ell-1}$ are pairwise distinct, all consecutive terms are adjacent in $G$,
and $v_\ell=v_j$ for some $0\leq j\leq\ell-3$.
We call $\ell=\ell(\beta)$ the length of the sequence.
Each such sequence corresponds to an extension excluded from the path-tree construction because it visits a previously visited vertex
(extensions that immediately return to the parent are not included).

\begin{lemma}\label{lem:root-defect}
Let $G$ be an $r$-regular graph on $n$ vertices $(r\ge3)$,
and let $x>2\sqrt{2(r-1)}$.
Setting $z=r$ in \eqref{eq:yqg}--\eqref{eq:tree-expression}
and using the above notation, we have
\begin{equation}\label{eq:root-defect}
 0\leq g-g_v\leq y g^2\sum_{\beta\in\mathcal{B}_v}q^{\ell(\beta)}.
\end{equation}
and
\begin{equation}\label{eq:blocked-count}
 \sum_{v\in V(G)}\sum_{\beta\in\mathcal{B}_v}q^{\ell(\beta)}
 \leq\sum_{C\in\mathcal{C}(G)}
       2|C|q^{|C|}\frac{1-q}{1-(r-1)q}.
\end{equation}
\end{lemma}

\begin{proof}
For a nonroot vertex $\omega$, let $\delta_\omega=h-h_\omega\geq0$.
By \eqref{eq:yqg}, \eqref{eq:y-identities}, and \eqref{eq:h-recursion},
\begin{align}
 \delta_\omega
 &=\frac1{x-2(r-1)h} - \frac1{x-2\sum_{\omega'\text{ is }\omega\text{'s child}}h_{\omega'}} \notag \\
 &=2hh_\omega\left((r-1)h-\sum_{\omega'}h_{\omega'}\right) \notag \\
 &=2hh_\omega\left(k_\omega h+a_\omega h-\sum_{\omega'}h_{\omega'}\right) \notag \\
 &=2hh_\omega\left(\sum_{\omega'\text{ is }\omega\text{'s child}}\delta_{\omega'} + a_\omega h\right) \notag \\
 &\leq q\sum_{\omega'\text{ is }\omega\text{'s child}}\delta_{\omega'} + q a_\omega h. \label{eq:defect-recursion}
\end{align}
Since the path tree is finite, \eqref{eq:defect-recursion} can be expanded all the way to the leaves.

For a nonroot node $\omega$ in the path tree,
its depth is its distance from the root $(v)$.
For a nonroot node at depth $k$,
$\omega=(v_0,\ldots,v_k)$,
the $a_\omega$ missing children correspond to vertices $v_j$ satisfying
$v_jv_k\in E(G)$ and $0\le j\le k-2$.
Thus each missing child uniquely produces
$\beta=(v_0,\ldots,v_k,v_j)\in \mathcal{B}_v$,
with $\ell(\beta)=k+1$.
In \eqref{eq:defect-recursion}, this missing child initially contributes $qh$;
propagating from the node at depth $k$ to a child of the root requires $k-1$ further applications of the recurrence,
each multiplying the contribution by at most $q$. Hence the total contribution is at most
$hq^k=hq^{\ell(\beta)-1}$.
Since $\mathcal{T}_v$ is a finite tree, summing over all missing children gives
\[
 \sum_{\omega\text{ is a child of the root}}\delta_\omega
 \leq h\sum_{\beta\in\mathcal{B}_v}q^{\ell(\beta)-1}.
\]
Since the root has exactly $r$ children and $2h/q=y$, we have
\begin{align*}
 g-g_v
 =2gg_v\sum_{\omega\text{ is a child of the root}}\delta_\omega
 \leq2g^2h\sum_{\beta\in\mathcal{B}_v}q^{\ell(\beta)-1}
 =yg^2\sum_{\beta\in\mathcal{B}_v}q^{\ell(\beta)}.
\end{align*}
This proves \eqref{eq:root-defect}.

We next count the sequences in \eqref{eq:blocked-definition}.
For a given $\beta$, the position $j$ of the repeated vertex is uniquely determined.
The part of the sequence from $v_j$ to $v_\ell=v_j$ is a cycle of length $L=\ell-j\geq3$,
while the part from $v_0$ to $v_j$ is a simple path of length $p=j$.
This path is vertex-disjoint from the cycle except at its endpoint $v_j$.
We call $v_j$ the attachment vertex of $\beta$ with respect to the cycle $C$.

For a fixed cycle $C$, there are $L$ choices of the attachment vertex;
there are two choices of traversal direction along $C$ from the attachment vertex.
When $p=0$, no attached path is selected.
When $p\geq1$, construct the attached path backwards from the attachment vertex (with at most $r-2$ choices for the first step
and at most $r-1$ choices for each subsequent step).
The length of an admissible attached path is at most $n-L$. Since all terms are nonnegative, for convenience in estimating the sum,
we may extend the finite sum to an infinite sum over $p\ge1$.
Thus, summing over all roots, the contribution of this cycle is at most
\begin{equation}\label{eq:jishu}
 2Lq^L\left(1+\sum_{p\geq1}(r-2)(r-1)^{p-1}q^p\right).
\end{equation}
Since $q<1/(r-1)$, the infinite series in \eqref{eq:jishu} is
a geometric series with common ratio $(r-1)q$, and it converges. Therefore
\[
 2Lq^L\left(1+\sum_{p\geq1}(r-2)(r-1)^{p-1}q^p\right)
 =2Lq^L\frac{1-q}{1-(r-1)q}.
\]
Summing over all cycles proves \eqref{eq:blocked-count}.
\end{proof}

\begin{lemma}[\cite{va3}]\label{lemma:30}
Let $G$ be a graph on $n$ vertices,
and let $m_k(G)$ denote the number of matchings of $G$ containing $k$ edges.
Define the matching polynomial by
\[
\mu(G,x)
=
\sum_{k=0}^{\lfloor n/2\rfloor}
(-1)^k m_k(G)x^{\, n-2k}.
\]
Then
\[
\frac{d}{dx}\mu(G,x)
=
\sum_{v\in V(G)}\mu(G-v,x).
\]
\end{lemma}

\begin{lemma}\label{lemma:matching-excess}
Let $G$ be an $r$-regular graph on $n$ vertices $(r\geq3)$. If $x>2\sqrt{2(r-1)}$,
then, with the notation in \eqref{eq:yqg} and \eqref{eq:tree-expression}, the following inequalities hold.
\begin{align}
 0\leq\ln M_G(x)-n\ln\mathcal{A}_r(x)
 \leq\sum_{C\in\mathcal{C}(G)}
    |C|\int_0^{q(x)}\frac{t^{|C|-1}}{1-t}\, dt
 \leq\frac1{1-q(x)}\sum_{C\in\mathcal{C}(G)}q(x)^{|C|}. \label{eq:matching-excess}
\end{align}
\end{lemma}
\begin{proof}
Define
\[
g_{v_i}(s)=\frac{M_{G-v_i}(s)}{M_G(s)},
\qquad
g(s)=\frac{y(s)}{y(s)^2-2}.
\]
By Lemma \ref{lemma:30},
$M_G'(x)=\sum_{v_i\in V(G)}M_{G-v_i}(x)$.
Combining this with \eqref{eq:root-ratio}, we have
\begin{equation}\label{eq:matching-log-derivative}
\frac{\mathrm{d}}{\mathrm{d}x} \ln M_G(x) = \sum_{v_i \in V(G)} g_{v_i}(x).
\end{equation}

On the other hand, differentiating $x$ and $\ln\mathcal{A}_r(x)$ with respect to $y$ gives
\[
 \frac{\mathrm{d}}{\mathrm{d}y}\ln\mathcal{A}_r(x(y))
 =\frac{y^2-2(r-1)}{y(y^2-2)},
 \qquad
 \frac{\mathrm{d} x}{\mathrm{d} y}=1-\frac{2(r-1)}{y^2},
\]
The condition $x>2\sqrt{2(r-1)}$ implies $y^2>2(r-1)$, so $dx/dy>0$.
Thus the chain rule gives
\begin{equation}\label{eq:tree-log-derivative}
 \frac{\mathrm{d}}{\mathrm{d}x}\ln\mathcal{A}_r(x)=\frac{y}{y^2-2}=g(x).
\end{equation}
As $x\to\infty$, we have $M_G(x)/x^n\to1$ and $\mathcal{A}_r(x)/x\to1$.
Hence $\ln M_G(x)-n\ln\mathcal{A}_r(x)\to0$ as $x\to\infty$.
Therefore, by \eqref{eq:matching-log-derivative}, \eqref{eq:tree-log-derivative}, and $g_{v_i}\leq g$,
\[
 \ln M_G(x)-n\ln\mathcal{A}_r(x)
 =\int_x^\infty\sum_{v_i\in V(G)}\bigl(g(s)-g_{v_i}(s)\bigr)\, ds\geq0.
\]
For every $s\ge x$, we have $s>2\sqrt{2(r-1)}$.
Applying Lemma \ref{lem:root-defect} at the parameter $s$ gives
\[
\int_x^\infty\sum_{v_i\in V(G)}\bigl(g(s)-g_{v_i}(s)\bigr)\,\mathrm{d} s \leq
\sum_{C\in\mathcal{C}(G)} \int_x^\infty y(s)g(s)^2
 \frac{2|C|q(s)^{|C|}(1-q(s))}{1-(r-1)q(s)}\, \mathrm{d}s.
\]
Now make the change of variables $u=y(s)$, $q(s)=2/u^2$.
Thus
\[
 \mathrm{d} s=(1-(r-1)q(s))\,\mathrm{d} u,
 \qquad y(s)g(s)^2=\frac1{u(1-q(s))^2}.
\]
Since $0\leq t\leq q(x)<1$, we have
$\frac1{1-t}\le\frac1{1-q(x)}$.
Therefore
\[
|C|
\int_0^{q(x)}
\frac{t^{|C|-1}}{1-t}\, dt
\le
\frac{|C|}{1-q(x)}
\int_0^{q(x)}
t^{|C|-1}\, dt
=
\frac{q(x)^{|C|}}{1-q(x)}.
\]
Hence
\begin{align*}
&\int_x^\infty y(s)g(s)^2
 \frac{2|C|q(s)^{|C|}(1-q(s))}{1-(r-1)q(s)}\, \mathrm{d}s \notag \\
&=\int_{y(x)}^\infty
       \frac{2|C|(2/u^2)^{|C|}}{u(1-2/u^2)}\,\mathrm{d} u \notag \\
&=|C|\int_0^{q(x)}\frac{t^{|C|-1}}{1-t}\,\mathrm{d} t \notag \\
&\leq\frac{q(x)^{|C|}}{1-q(x)}. 
\end{align*}
Summing over the finite set $\mathcal{C}(G)$ yields \eqref{eq:matching-excess}.
\end{proof}

\section{Proof of Theorem \ref{thm2}}
In this section, we give the proof of Theorem \ref{thm2} as follows.

\begin{proofthm}{Proof of Theorem \ref{thm2}}
For $r=3$, the result was already proved by Bencs and Csikvári \cite{Bencs2023}.

Henceforth, assume $r\geq4$ and set $x_0=r+1$.
Since $x_0^2-8(r-1)=(r-3)^2>0$,
Lemma \ref{lemma:matching-excess} applies, and
\[
 y(x_0)=r-1, \qquad q(x_0)=q_r=\frac2{(r-1)^2},
 \qquad \mathcal{A}_r(x_0)=\frac{(r - 1)^{r-1}}{(r^2 - 2r - 1)^{r/2-1}}.
\]
Combining \eqref{eq:cycle-loss} with \eqref{eq:matching-excess} gives
\begin{equation}\label{eq:combine}
 \ln F(G)-n\ln \mathcal{A}_r(x_0)
 \leq\sum_{C\in\mathcal{C}(G)}
       \left(\frac{q_r^{|C|}}{1-q_r}-2\rho_r^{-|C|}\right).
\end{equation}

We now compare the constants.
First, $0<q_r\leq2/9<1/2$.
Second,
\begin{align*}
 (r+1)(r-1)^2-2(r^2+r+1)
 &=r^3-3r^2-3r-1\\
 &=(r-4)^3+9(r-4)^2+21(r-4)+3>0.
\end{align*}
so
\[
 q_r\rho_r
 =\frac{2(r^2+r+1)}{(r+1)(r-1)^2}<1.
\]
Therefore, for every integer $b\geq3$,
\begin{align}
 \frac{q_r^b}{1-q_r}-2\rho_r^{-b}
 \leq\left(\frac1{1-q_r}-2\right)\rho_r^{-b}
 =-\eta_r \rho_r^{-b}, \label{eq:termwise-comparison}
\end{align}
where $\eta_r=(1-2q_r)/(1-q_r)>0$.
Substituting \eqref{eq:termwise-comparison} into \eqref{eq:combine} gives
\begin{equation}\label{eq:quantitative-log-bound}
 \ln F(G)\leq n\ln \mathcal{A}_r(x_0)
       -\eta_r\sum_{C\in\mathcal{C}(G)}\rho_r^{-|C|}.
\end{equation}

Since an $r$-regular graph $G$ with $r\geq4$ must contain a cycle, we have $\mathcal{C}(G)\neq\varnothing$.
Thus
\[
 \ln F(G)\leq n\ln \mathcal{A}_r(x_0)
       -\eta_r\sum_{C\in\mathcal{C}(G)}\rho_r^{-|C|}
       <n\ln \mathcal{A}_r(x_0)
       =n\ln \frac{(r - 1)^{r-1}}{(r^2 - 2r - 1)^{r/2-1}}.
\]
Therefore $F(G)^{1/n}<\frac{(r - 1)^{r-1}}{(r^2 - 2r - 1)^{r/2-1}}$,
completing the proof.
\end{proofthm}

The proof yields the following quantitative form of the upper bound.
\begin{corollary}\label{cor:cycle-correction}
Let $G$ be an $r$-regular graph on $n\geq1$ vertices, where $r\geq4$, and set
\[
q_r=\frac{2}{(r-1)^2},\qquad
\rho_r=r+1-\frac{r}{r+1},\qquad
\eta_r=\frac{1-2q_r}{1-q_r}>0.
\]
Then
\[
F(G)\leq
\left(\frac{(r-1)^{r-1}}{(r^2-2r-1)^{r/2-1}}\right)^n
\exp\!\left(-\eta_r\sum_{C\in\mathcal{C}(G)}\rho_r^{-|C|}\right).
\]
\end{corollary}
\begin{proof}
Exponentiate \eqref{eq:quantitative-log-bound} and use the displayed evaluation of $\mathcal{A}_r(x_0)$ in the proof of Theorem \ref{thm2}.
\end{proof}

\section{Conclusion}

This paper investigates the problem of counting spanning forests of 
r-regular graphs and establishes an upper bound by employing a weighted generalization of the matrix-forest theorem, graph polynomials, path trees, and related techniques. These results further elucidate the influence of vertex degrees and cycle structure on the number of spanning forests, and lay the groundwork for studying upper bounds on the number of spanning forests in more general graphs.

\noindent{\bf Conflict of Interest Statement}

The authors declare that they have no conflicts of interest.
\renewcommand\refname{\ References}

\end{document}